\documentclass{amsart}
\usepackage{amssymb}
\usepackage{amsthm}
\begin{document}
\newtheorem{defn}{Definition}[section]
\newtheorem{thm}{Theorem}[section]
\newtheorem{pr}{Proposition}[section]
\newtheorem{exam}{Example}[section]
\newtheorem{cor}{Corollary}[section]
\newtheorem{rem}{Remark}[section]
\newtheorem{lem}{Lemma}[section]
\newenvironment{dem}{\rm \trivlist \item[\hskip \labelsep{\bf
Proof}.]}{\par\nopagebreak \hfill $\Box$ \endtrivlist}
\date{}
\author{J.M. Cabezas, L.M. Camacho, J.R. G\'{o}mez, B.A. Omirov}
\title{\bf On the description of the Leibniz algebras with nilindex
$n-3$}
\maketitle
\begin{abstract}
In this paper we present the classification of a subclass of
naturally graded Leibniz algebras. These $n$-dimensional Leibniz
algebras have the characteristic sequence
equal to $(n-3,3).$ For this purpose we use the software
$Mathematica$.
\end{abstract}
\medskip \textbf{AMS Subject Classifications (2000):
17A32, 17A36, 17A60, 17B70.}
\textbf{Key words:} Lie algebra, Leibniz algebra, nilpotence,
natural gradation, cha\-rac\-te\-ris\-tic sequence,
$p$-filiformlicity.
\def\gg{{\mathfrak g}}
\def\CC{{\mathbb{C}}}
\def\ll{{ L}}
\section{Introduction}
Leibniz algebras are one of the new algebras introduced by Loday
\cite{loday}, \cite{loday2} in connection with the study of
periodicity phenomena in algebraic K-theory. Leibniz algebras
have been introduced as a "non-antisymmetric" analogue of Lie
algebras.
A Leibniz algebra $L$ is a vector space equipped with a bracket
[-,-] satisfying the identity $$[x, [y, z]] = [[x, y], z]-[[x, z],
y].$$

If the antisymmetric relation is assumed, this identity is
equivalent to the Jacobi identity. Hence, a Lie algebra is a
Leibniz algebra.
It is well known that the natural gradation of nilpotent Lie and
Leibniz algebras is very helpful in investigating their
structural properties. A remarkable fact of the naturally graded
algebras is the relative simplicity of the study of the
cohomological properties, (see for example \cite{Dz1}- \cite{G-H}
and \cite{ Ve}).

Recently, some papers are focused to the study of some interesting
families of Leibniz algebras, such as $p$-filiform and
quasi-filiform Leibniz algebras. These algebras have their
characteristic sequences equal to $(n-p,1,1,...,1)$ and $(n-2,2)$
with $dim(\ll)=n$, \cite{Comm-p-fil}--\cite{JSC}.

Naturally graded $p$-filiform Leibniz algebras are already
classified in \cite{J.Lie.Theory} and \cite{Comm-p-fil}. The
classification of naturally graded nul-filiform and filiform
Leibniz algebras reader can find in \cite{Omirov1}. The
quasi-filiform $n$-dimensional Leibniz algebras have
characteristic sequence $(n-2,1,1)$ (the case of $2$-filiform) or
\mbox{$(n-2,2)$} \cite{Comm-2-fil} and \cite{JSC}.

For a given Leibniz algebra $\ll$ we define the descending central
series as follows:
$$L^1=L, \quad L^{k+1}=[L^k,L], \quad k\geq 1.$$

If there exists a natural number $s$ such that $L^s={0},$ then the
Leibniz algebra $L$ is said to be nilpotent and minimal such
number is called the nilindex of the algebra $L.$

Bellow we present a gradation closely related to the descending
central series.

Let $L$ be a nilpotent Leibniz algebra with nilindex $s.$ We put
$L_i=L^i/L^{i+1}$ for $1\leq i\leq s-1,$ and $gr L=L_1\oplus
L_2\oplus \dots \oplus L_{s-1}.$ It is easy to check embedding
$[L_i,L_j]\subseteq L_{i+j}$ and therefore, the algebra $gr L$ is
graded algebra, which is called the naturally graded Leibniz
algebra.

Let $x$ be a nilpotent element of the set $L \setminus L^2$. For
the nilpotent operator of right multiplication $R_x$ we define a
decreasing sequence $C(x)=(n_1,n_2, \dots, n_k)$, which consists
of the dimensions of Jordan blocks of the operator $R_x$. On the
set of such sequences we consider the lexicographic order, that
is,
$C(x)=(n_1,n_2, \dots, n_k)\leq C(y)=(m_1, m_2, \dots, m_s)\Longleftrightarrow$
there exists $i \in \mathbb{N}$ such that $n_j=m_j$ for any $j<i$ and $n_i<m_i$.

The sequence $C(\ll)=max \ C(x)_{x \in \ll \setminus \ll^2}$ is
called characteristic sequence of the algebra $\ll$.
If $C(\ll)=(1,1, \dots,1)$ then evidently, the algebra $\ll$ is abelian.

The set $R(L)=\{x\in L \ | \ [y,x]=0 \ \mbox{for \ any} \ y\in
L\}$ is said to be a right annihilator of the algebra $L.$

In this work we classify a subclass of naturally graded
Leibniz algebras with nilindex $n-3.$ In case of Leibniz algebras
with nilindex equal to $n-3,$ for the characteristic sequence we
have the following tree possibilities: $$(n-3,1,1,1), \ (n-3,2,1)
\ \mbox{and} \ (n-3,3).$$
The first one is $3$-filiform case. We
will focus our attention on the study of those with characteristic
sequence $(n-3, 3)$. Throughout all the work, we use the software
$Mathematica$. Since in the case of non-Lie Leibniz algebras the
skew-symmetric identity is not valid, this classification is
very complex and we should overcome the difficulties, which
need a lot of computations. Using computer programs is very
helpful for computing the Leibniz identity in low dimension and
formulate the generalizations of the calculations, which are
proved for arbitrary finite dimension. The used program can be
find in \cite{JSC}. Some examples of the programs for various
types of Leibniz algebras classes are in the following Web site:
http://personal.us.es/jrgomez.
\section{Naturally graded
Leibniz algebras with characteristic sequence $(n-3,3)$.}
Let $\ll$ be a naturally graded $n$-dimensional Leibniz algebra
which characteristic sequence equal to $(n-3,3)$. From the
definition of the characteristic sequence, it follows the
existence of a basis $\{e_{1}, e_{2},\dots, e_{n}\}$ such that
element $e_{1}\in \ll\backslash \ll^2$ and the operator of right
multiplication $R_{e_1}$
has one of the following
forms: $$\begin{array}{ll} \left(\begin{array}{ll}
J_{n-3}&0\\
0&J_{3}
\end{array}\right),&
\left(\begin{array}{ll}
J_{3}&0\\
0&J_{n-3}
\end{array}\right)
\end{array}$$
\begin{defn}
A naturally graded Leibniz algebra $\ll$ which characteristic
sequence is equal to $(n-3,3)$, is called algebra of the second
type if there exists a basic element $e_{1}\in \ll\backslash
\ll^2$ such that the operator $R_{e_1}$ has the form:
$$\left(\begin{array}{ll}
J_{n-3}&0\\
0&J_{3}
\end{array}\right);$$
if $R_{e_1}$ has the other form, then it is called algebra of the
second type.
\end{defn}
Since the classification of Leibniz algebras of the second type is
more complicated and it needs to use more original technics, first we
present the description of the second type.
\begin{thm}
Let $\ll$ be an $n$-dimensional naturally graded Leibniz algebra
of the second type ($n\geq 9$). Then it is isomorphic to one
of the following pairwise non-isomorphic algebras:
$$\small\begin{array}{|l|l|l|l|}
\hline
&\lambda&\mu& dim(L)\\[1mm]
\hline
\hline
\ll^{0,1}_{(0,0,0,0,0)}& & & \mbox{odd or even}\\[1mm]
\hline
\ll^{0,2}_{(0,0,0,\lambda,-1)}&\lambda \in \{0,1\} & & \mbox{odd or even}\\[1mm]
\hline
\ll^{0,3}_{(1,0,0,\lambda,-1)}&\lambda \in \mathbb{C} & & \mbox{odd or even}\\[1mm]
\hline
\ll_{(1,0,1/4,\lambda,-1)}^{0,4}&\lambda \in \mathbb{C} & & \mbox{odd or even}\\[1mm]
\hline
\ll_{(0,0,1,\lambda,-1)}^{0,5}&\lambda \in \mathbb{C}& & \mbox{odd or even}\\[1mm]
\hline
\ll_{(0,1,0,\lambda,-1)}^{0,6}&\lambda \in \{0,1\}& & \mbox{odd or even}\\[1mm]
\hline
\ll_{(\mu,1,0,\lambda,-1)}^{0,6}&\lambda \in \mathbb{C}&\mu \in \{1,2\} & \mbox{odd or even}\\[1mm]
\hline
\ll_{(0,1,\mu,\lambda,-1)}^{0,7}&\lambda \in \mathbb{C}&\mu \in \mathbb{C}\setminus\{0\}
& \mbox{odd or even}\\[1mm]
\hline
\ll_{(-2\lambda,1,-\lambda,2,-1)}^{0,8}&\lambda \in \{-2,-4/3\} & & \mbox{odd or even}\\[1mm]
\hline
\ll_{(2\lambda,1,\lambda,0,-1)}^{0,9}&\lambda \in \mathbb{C}\setminus\{0,1\}& &
\mbox{odd or even}\\[1mm]
\hline
\ll_{(1,1,1/4,1/4,-1)}^{0,10}& & & \mbox{odd or even}\\[1mm]
\hline
\ll_{(1,1,1/4,1/2,-1)}^{0,10}& & & \mbox{odd or even}\\[1mm]
\hline
\ll_{(2,1,1,1,-1)}^{0,10}& & & \mbox{odd or even}\\[1mm]
\hline
\ll_{(2,1,1,0,-1)}^{0,10}& & & \mbox{odd or even}\\[1mm]
\hline
\ll_{(1,\lambda,1/4,0,-1)}^{0,11}&\lambda \in \mathbb{C}\setminus\{0,1/2\}& & \mbox{odd or even}\\[1mm]
\hline
\ll^{1,2}_{(0,0,0,\lambda,-1)}&\lambda \in \{0,1\} & & \mbox{even}\\[1mm]
\hline
\ll^{1,3}_{(1,0,0,\lambda,-1)}&\lambda \in \mathbb{C} & & \mbox{even}\\[1mm]
\hline
\ll_{(1,0,1/4,\lambda,-1)}^{1,4}&\lambda \in \mathbb{C} & & \mbox{even}\\[1mm]
\hline
\ll_{(\mu,1,0,\lambda,-1)}^{1,6}&\lambda \in \mathbb{C}&\mu \in \mathbb{C}& \mbox{even}\\[1mm]
\hline
\ll_{(0,\gamma,\mu,\lambda,-1)}^{1,7}&\lambda \in \mathbb{C}&\gamma,\
\mu \in \mathbb{C}\setminus\{0\} & \mbox{even}\\[1mm]
\hline
\ll_{(-2\lambda,1,\lambda,\mu,-1)}^{1,9}&\lambda \in \mathbb{C}\setminus\{0,1\} & \mu \in \mathbb{C}& \mbox{even}\\[1mm]
\hline
\ll_{(\lambda,1,\lambda^2/4,\mu,-1)}^{1,11}&\lambda \in
\mathbb{C}\setminus\{-2,0\} &\mu \in \mathbb{C} & \mbox{even}\\[1mm]
\hline
\ll_{(-1,0,0,\lambda,-1)}^{1,12}&\lambda \in \{0,1\} & & \mbox{even}\\[1mm]
\hline
\ll_{(-2,0,1,\lambda,-1)}^{1,13}&\lambda \in \CC & & \mbox{even}\\[1mm]
\hline
\ll_{(-4,0,2,\lambda,-1)}^{1,14}&\lambda \in \CC & & \mbox{even}\\[1mm]
\hline
\ll_{(0,0,-1,\lambda,-1)}^{1,15}&\lambda \in \CC & & \mbox{even}\\[1mm]
\hline
\ll_{(-2,0,-1,\lambda,-1)}^{1,16}&\lambda \in \CC & & \mbox{even}\\[1mm]
\hline
\ll_{(0,-1,0,\lambda,-1)}^{1,17}&\lambda \in \{0,1\} & & \mbox{even}\\[1mm]
\hline
\ll_{(-1,-1,0,\lambda,-1)}^{1,18}&\lambda \in \CC & & \mbox{even}\\[1mm]
\hline
\ll_{(-2,-1,0,1,-1)}^{1,19}& & & \mbox{even}\\[1mm]
\hline
\ll_{(1,-1,0,\lambda,-1)}^{1,20}&\lambda \in \CC\setminus\{-1/2\} & & \mbox{even}\\[1mm]
\hline
\ll_{(1,1/3,0,\lambda,-1)}^{1,21}&\lambda \in \CC & & \mbox{even}\\[1mm]
\hline
\ll_{(-2,-1,1,\lambda,-1)}^{1,22}&\lambda \in \{0,1\} & & \mbox{even}\\[1mm]
\hline
\ll_{(1,1/2,1/4,\lambda,-1)}^{1,23}&\lambda \in \CC & & \mbox{even}\\[1mm]
\hline
\end{array}$$
$$\small\begin{array}{|l|l|l|c|}
\hline
&\lambda&\gamma, \ \mu& dim(L)\\[1mm]
\hline
\hline
\ll_{(-4,-1,2,\lambda,-1)}^{1,24}&\lambda \in \CC & & \mbox{even}\\[1mm]
\hline
\ll_{(-3,-4/3,2,\lambda,-1)}^{1,25}&\lambda \in \CC & & \mbox{even}\\[1mm]
\hline
\ll_{(2/5,2,2/5,\lambda,-1)}^{1,26}&\lambda \in \CC & & \mbox{even}\\[1mm]
\hline
\ll_{(2/\lambda,\lambda,1,\mu,-1)}^{1,27}&\lambda \in
\CC\setminus\{-1,0,1\} &\mu\in \CC & \mbox{even}\\[1mm]
\hline
\ll_{(8/5,1/2,-4/5,\lambda,-1)}^{1,28}&\lambda \in \CC & & \mbox{even}\\[1mm]
\hline
\ll_{(\lambda,-1,\lambda^2/4,0,-1)}^{1,29}&\lambda \in \CC\setminus\{-2,0\} & & \mbox{even}\\[1mm]
\hline
\ll_{(1,-1,1/4,\lambda,-1)}^{1,30}&\lambda \in \{-1/2,1/4\} & & \mbox{even}\\[1mm]
\hline
\ll_{(-8,2,16,\lambda,-1)}^{1,31}&\lambda \in \CC & & \mbox{even}\\[1mm]
\hline
\ll_{(-2,\lambda,1,0,-1)}^{1,32}&\lambda \in \CC\setminus\{-1,0\} & & \mbox{even}\\[1mm]
\hline
\ll_{(-2,1,1,\lambda,-1)}^{1,33}&\lambda \in \{-1,1\} & & \mbox{even}\\[1mm]
\hline
\end{array}$$
where the algebra
$$\small\begin{array}{l}
\ll^{\epsilon,j}_{(\alpha_1,\alpha_2,\alpha_3,\alpha_4,\beta)}:\quad \epsilon\in \{0,1\}, \quad 1\leq j\leq 33,\quad \beta\in\{-1,0\}\\
\end{array}$$
has the following multiplication:
$$\left\{\begin{array}{l}
\left[e_i,e_1\right]=e_{i+1},\ 1\leq i\leq n-1,\ i\neq 3\\
\left[e_1,e_4\right]=\alpha_1 e_2+\beta e_5,\\
\left[e_2,e_4\right]=\alpha_2 e_3,\\
\left[e_4,e_4\right]=\alpha_3 e_2,\\
\left[e_5,e_4\right]=\alpha_4 e_3,\\
\left[e_1,e_5\right]=(\alpha_1-\alpha_2) e_3-e_6,\\
\left[e_4,e_5\right]=(\alpha_3-\alpha_4) e_3,\\
\left[e_1,e_i\right]=\beta e_{i+1},\ 6\leq i \leq n-1,\\
\left[e_i,e_{n+3-i}\right]=\epsilon (-1)^ie_n,\ 4 \leq i \leq n-1.
\end{array}\right.
$$
\end{thm}
\dem
From the condition of the theorem we have the
following multiplication of the basic element $e_1$ on the right
side: $$[e_i,e_1]=e_{i+1},\ 1\leq i\leq n-1,\ i\neq 3, \
[e_3,e_1]=[e_n,e_1]=0.$$

From these products we conclude that $$\ll_1=<e_1, e_4>, \
\ll_2=<e_2, e_5>, \ \ll_3=<e_3, e_6>, \ \ll_i=<e_{i+3}>, \
4\leq i\leq n-3$$
and $e_2, e_3\in R(\ll).$

Let us introduce denotations
$$\small\begin{array}{lll} [e_1,e_4]=\alpha_1 e_2+\beta_1
e_5,&[e_2,e_4]=\alpha_2 e_3+\beta_2 e_6,\quad &[e_3,e_4]=\beta_3
e_7,\\{}
[e_4,e_4]=\alpha_3 e_2+\beta_4
e_5,&[e_5,e_4]=\alpha_4 e_3+\beta_5 e_6,\quad &\\{}
[e_i,e_4]=\beta_i e_{i+1},\ 6\leq i\leq n-1,&[e_n,e_4]=0.&
\end{array}$$

The equalities $[e_i,e_5]=[[e_i,e_4],e_1]-[[e_i,e_1],e_4],\ 1\leq
i\leq n$ derive
$$\small\begin{array}{lll} [e_1,e_5]=(\alpha_1-\alpha_2)e_3+(\beta_1-\beta_2)e_6,&[e_2,e_5]=(\beta_2-\beta_3)e_7, &
[e_3,e_5]=\beta_3 e_8, \\{}
[e_4,e_5]=(\alpha_3-\alpha_4)e_3+(\beta_4-\beta_5)e_6,&[e_5,e_5]=(\beta_5-\beta_6)e_7,
& \\{} [e_i,e_5]=(\beta_i-\beta_{i+1})e_{i+2},\ 6\leq i\leq n-2&
[e_{n-1},e_5]=[e_n,e_5]=0.&
\end{array}$$

Using induction on $j$ for any value $i$ it can be proved that
$$ [e_{i},e_{j}]=\displaystyle\left(\sum_{k=0}^{j-4}(-1)^k
\left(\begin{array}{c}
j-4\\
k
\end{array}\right) \beta_{i+k}\right)e_{i+j-3},\ \ 5 \leq i \leq n-3,\ \ 6 \leq j \leq n+3-i. $$

In the case of $e_4\in R(\ll)$ we obtain the algebra
$\ll^{0,1}_{(0,0,0,0,0)}.$

Let now $e_{4}\notin R(\ll).$ Then we consider the following cases:

\fbox{$e_{5}\in R(\ll)$}

Then $e_i\in R(\ll)$ for $2\leq i \leq
n, \ i\neq 4$.

From the equalities $[[e_i,e_1],e_4]=[[e_i,e_4],e_1],\ 1\leq i
\leq n$, we have $$\alpha_2=\alpha_1,\ \alpha_4=\alpha_3,\
\beta_3=\beta_2=\beta_1,\ \beta_i=\beta_4, \  5\leq i \leq
n-1.$$
For $n\geq 8$ we have also $\beta_1=0.$

The change of basis taken as
$$e'_i=e_i,\ 1\leq i \leq n, \ i\neq
4,5,6 ,\quad e'_{j}=e_j-\beta_4e_{j-3},\ 4\leq j \leq 6$$
deduces $\beta_4=0$.

If we take the change of basis in the following way:
$$e'_1=Ae_1+Be_4,\ e'_{n-2}=e_1,\ e'_j=[e'_{j-1},e'_1],\ 2\leq
j\leq n, \ j\neq n-2$$ with condition $AB(A+\alpha_1B)\neq 0$,
then we obtain the algebra of the first type. Therefore, this case
is impossible for the algebra of the second type.

\fbox{$e_{5}\notin R(\ll)$}

The embedding $[e_4,e_4]\in R(\ll)$ implies $\beta_4=0$ and from
$[e_i,[e_4,e_1]]=-[e_i,[e_1,e_4]],$ with $1\leq i \leq n$ we obtain$\beta_1=-1$.

If $e_6\in R(\ll)$, then for $n\geq 9$ it follows $\beta_1=0,$
which is a contradiction with the condition $\beta_1=-1$.
Therefore, $e_6\notin R(\ll).$

It is easy to check that $[e_i,e_j]+[e_j,e_i]\in R(\ll)$ for any
values of $i,j$. Applying this for $i=1$ and $j=5$ we
obtain $\beta_2=0$.

The following equalities:
$$[e_1,e_i]=-e_{i+1},\quad [e_2,e_i]=[e_3,e_i]=0,\ 6\leq i \leq
n-1$$ are proved by induction on $i.$

From $[e_1,[e_4,e_{2j+1}]]=-[e_5,e_{2j+1}]+[e_{2j+2},e_4],\ j\geq
2 $, we have that
$$2\beta_{2j+2}=\beta_5+\beta_{2j+1}+\sum_{k=1}^{2j-4}(-1)^k
\left(\begin{array}{c}
2j-3\\
k
\end{array}\right) (\beta_{5+k}-\beta_{4+k}), \ j\geq 2.$$

Similar as in \cite{JSC} we derive
$$\left\{\begin{array}{ll}
\beta_j=\beta_5,&6\leq j\leq n-1, \ \mbox{for n odd,}\\
\beta_j=\beta_5,& 6\leq j\leq n-2, \ \mbox{for n even}
\end{array}\right.$$
and


\hspace{2cm}$[e_4,e_{n-1}]=-\beta_5 e_n$ for $n$ odd,

\hspace{2cm}$[e_4,e_{n-1}]=(\beta_{n-1}-2\beta_5)e_n$ for $n$ even,

\hspace{2cm}$[e_i,e_{n+3-i}]=(-1)^i (\beta_{n-1}-\beta_5)e_n,\
5\leq i \leq n-2$, for $n$ even.

\

If $\beta_{n-1}=\beta_5$, then by the change of basis defined as $
e'_i=e_i,\ 1\leq i \leq n, $ $ i\neq 4,5,6,$ and $e'_i=e_i-\beta_5
e_{i-3},$ $ 4 \leq i \leq 6$ we can assume $\beta_5=0$.

If $\beta_{n-1}\neq \beta_5$ (the case of $n$ even), then by using
the change of basis:
$$\left\{\begin{array}{ll}
e'_i=(\beta_{n-1}-\beta_5)^i e_i,& 1\leq i \leq 3,\\
e'_4=e_4-\beta_5 e_1,&\\
e'_5=(\beta_{n-1}-\beta_5)(e_{5}-\beta_5 e_2),&\\
e'_6=(\beta_{n-1}-\beta_5)^2(e_6-\beta_5 e_3),&\\
e'_i=(\beta_{n-1}-\beta_5)^{i-4}e_i,&7\leq i \leq n
\end{array}\right.$$
we obtain $\left[e_i,e_{n+3-i}\right]=(-1)^ie_n$ for $ 4 \leq i
\leq n-1.$
Thus, multiplication in $\ll$ is as follows:
$$\left\{\begin{array}{ll}
\left[e_i,e_1\right]=e_{i+1},& 1\leq i\leq n-1,\ i\neq 3,\\
\left[e_1,e_4\right]=\alpha_1 e_2-e_5,&\\
\left[e_2,e_4\right]=\alpha_2 e_3,&\\
\left[e_4,e_4\right]=\alpha_3 e_2,&\\
\left[e_5,e_4\right]=\alpha_4 e_3,&\\
\left[e_1,e_5\right]=(\alpha_1-\alpha_2) e_3-e_6,&\\
\left[e_4,e_5\right]=(\alpha_3-\alpha_4) e_3,&\\
\left[e_1,e_i\right]=-e_{i+1},& 6\leq i \leq n-1,\\
\left[e_i,e_{n+3-i}\right]=\epsilon (-1)^ie_n,& 4 \leq i \leq
n-1,\ \epsilon\in\{0,1\}.
\end{array}\right.$$
\
\fbox{{\bf Case 1.} $\epsilon=0$} ($n$ odd or even)

Applying the general change of generators of the basis:
\begin{center}
$e'_1=\displaystyle\sum_{i=1}^{n} A_i e_i,\qquad
e'_{n-2}=\sum_{i=1}^{n} B_i e_i, $
\end{center}
\noindent
we determine the other elements of the new basis and the products
in this basis. Then the new parameters are the following:
$$\small\begin{array}{l} \alpha'_1=\displaystyle\frac{(\alpha_1
A_1+2\alpha_3 A_4)B_4 }{A_1^2+\alpha_1 A_1A_4+\alpha_3
A_4^2},\quad \ \ \
\alpha'_2=\displaystyle\frac{\alpha_2B_4}{A_1+\alpha_2A_4},\\[5mm]
\alpha'_3=\displaystyle\frac{\alpha_3B_4^2}{A_1^2+\alpha_1
A_1A_4+\alpha_3 A_4^2}, \quad
\alpha'_4=\displaystyle\frac{(\alpha_4 A_1+\alpha_2 \alpha_3
A_4)B_4^2 }{(A_1+\alpha_2A_4)(A_1^2+\alpha_1 A_1A_4+\alpha_3
A_4^2)},\end{array}$$ satisfying the restriction
$A_1(A_1+\alpha_2A_4)(A_1^2+\alpha_1 A_1A_4+\alpha_3
A_4^2)B_4\neq 0$.

Note that for new parameters we have $$\begin{array}{l}
\alpha'^2_{1}-4\alpha'_3=\displaystyle\frac{(\alpha_1^2-4\alpha_3)A_1^2B_4^2}{(A_1^2+\alpha_1
A_1A_4+\alpha_3 A_4^2)^2}, \\[5mm]
\alpha'_1\alpha'_2-2\alpha'_3=\displaystyle\frac{(\alpha_1\alpha_2-2\alpha_3)A_1B_4^2}{(A_1+\alpha_2A_4)(A_1^2+\alpha_1
A_1A_4+\alpha_3 A_4^2)}, \\[5mm]
\alpha'_1\alpha'_2-2\alpha'_4=\displaystyle\frac{(\alpha_1\alpha_2-2\alpha_4)A_1B_4^2}{(A_1+\alpha_2A_4)(A_1^2+\alpha_1
A_1A_4+\alpha_3 A_4^2)}.
\end{array}$$
Consequently, the nullity of $\alpha_1^2-4\alpha_3$ is invariant
in the following sense:

if $\alpha_1^2-4\alpha_3=0,$ then ${\alpha'}_1^2-4{\alpha'}_3=0$
and if $\alpha_1^2-4\alpha_3\neq 0,$ then
${\alpha'}_1^2-4{\alpha'}_3\neq 0$. Analogously, the expressions
$\alpha_1\alpha_2-2\alpha_3 $ and $ \alpha_1\alpha_2-2\alpha_4$
are nullity invariants.

Consider the following subcases:

\noindent\fbox{\bf $\alpha_2=0$,\ $\alpha_3=0$ }

Then, $\alpha'_1=\displaystyle\frac{\alpha_1 B_4}{A_1+\alpha_1
A_4},$ $\alpha'_2=0,$ $\alpha'_3=0$ and
$\alpha'_4=\displaystyle\frac{ \alpha_4B_4^2}{A_1(A_1+\alpha_1
A_4)}.$
\begin{itemize}
\item $\alpha_1=0$.

If $\alpha_4=0,$ then the algebra $\ll^{0,2}_{(0,0,0,\lambda,-1)}$
with $\lambda=0$ is obtained.

If $\alpha_4\neq 0,$ then we obtain the algebra
$\ll^{0,2}_{(0,0,0,\lambda,-1)}$ with $\lambda=1$.
\item $\alpha_1\neq 0$.

If $\alpha_4=0,$ then we easily obtain $\alpha'_1=1$.
Thus, we have the algebra $\ll^{0,3}_{(1,0,0,\lambda,-1)}$ with $\lambda=0.$

If $\alpha_4\neq 0,$ then choosing appropriate values of $A_4$ and
$B_4$ we derive \mbox{$\alpha'_1=\alpha'_4=1$.} Hence, the algebra
$\ll^{0,3}_{(1,0,0,\lambda,-1)}$ with $\lambda=1$ is obtained.
\end{itemize}

\noindent\fbox{\bf $\alpha_2=0$,\ $\alpha_3\neq 0$}

Then, \quad $\small \alpha'_1=\displaystyle\frac{(\alpha_1
A_1+2\alpha_3 A_4)B_4 }{A_1^2+\alpha_1 A_1A_4+\alpha_3 A_4^2}, \
\ \ \ \alpha'_2=0,$

\
\hspace{1.4cm}$\small
\alpha'_3=\displaystyle\frac{\alpha_3B_4^2}{A_1^2+\alpha_1
A_1A_4+\alpha_3 A_4^2}, \ \ \ \
\alpha'_4=\displaystyle\frac{\alpha_4 B_4^2 }{A_1^2+\alpha_1
A_1A_4+\alpha_3 A_4^2} $.
\
\begin{itemize}
\item If $\alpha_1^2-4\alpha_3=0,$ then taking
adequate value of $B_4$ we obtain $\alpha'_1=1,$ $\alpha'_3=1/4$ and
$\alpha'_4=\displaystyle\frac{\alpha_4}{\alpha_1^2}=\lambda$. So,
we obtain the family of algebras
$\ll^{0,4}_{(1,0,1/4,\lambda,-1)}$ with $\lambda\in \CC.$
\item If $\alpha_1^2-4\alpha_3\neq0,$
then taking suitable values of $A_4$ and $B_4$
we deduce $\alpha'_1=0,$ $ \alpha'_3=1$ and
$\alpha'_4=\displaystyle\frac{\alpha_4}{\alpha_3}=\lambda$. The
family $\ll^{0,5}_{(0,0,1,\lambda,-1)}$, $\lambda\in \CC$ is
obtained.
\end{itemize}
\noindent\fbox{\bf $\alpha_2\neq 0$,\ $\alpha_3= 0$}

\

Then,
$$\small \alpha'_1=\displaystyle\frac{\alpha_1B_4 }{A_1+\alpha_1
A_4},\ \ \
\alpha'_2=\displaystyle\frac{\alpha_2B_4}{A_1+\alpha_2A_4},\ \
\alpha'_3=0,\ \alpha'_4=\displaystyle\frac{\alpha_4B_4^2
}{(A_1+\alpha_1A_4)(A_1+\alpha_2 A_4)}. $$
\begin{itemize}
\item $\alpha_1=0$.

If $\alpha_4=0,$ then the choosing appropriate $B_4$ leads
$\alpha'_2=1$. Thus, we obtain $\ll^{0,6}_{(0,1,0,\lambda,-1)},$ $
\lambda=0.$

If $\alpha_4\neq 0,$ then taking adequate $A_4$ and $B_4$ we
derive $\alpha'_2=\alpha'_4=1$. The algebra
$\ll^{0,6}_{(0,1,0,\lambda,-1)},$ $ \lambda=1$ is obtained.
\item $\alpha_1\neq 0$.
\begin{itemize}
\item[$\checkmark$] $\alpha_4=0$.

If $\alpha_1-\alpha_2= 0,$ then for suitable
$B_4$ we have $\alpha'_1=\alpha'_2=1,$ i.e.
we obtain the algebra $\ll^{0,6}_{(\mu,1,0,\lambda,-1)}$ with $\mu=1,\ \lambda=0.$

If $\alpha_1-\alpha_2\neq 0,$ then for adequate $A_4$ and
$B_4$ it follows that $\alpha'_1=2 $, $\alpha'_2=1$. The algebra
$\ll^{0,6}_{(\mu,1,0,\lambda,-1)},$ with $\mu=2, \lambda=0$ is
obtained.
\item[$\checkmark$] $\alpha_4\neq 0$.

If $\alpha_1-\alpha_2= 0,$ then for appropriate value of
$B_4$ we have $\alpha'_1=\alpha'_2=1$ and
$\alpha'_4=\displaystyle\frac{\alpha_4}{\alpha_1^2}=\lambda$.
Therefore, we obtain the family of algebras
$\ll^{0,6}_{(\mu,1,0,\lambda,-1)}$, where $\mu=1$, $\lambda\in
\CC\setminus\{0\}.$

If $\alpha_1-\alpha_2\neq 0,$ then taking suitable values of
$A_4$ and $B_4$ we obtain $\alpha'_1=2,$ $\alpha'_2=1$, \quad
$\alpha'_4=\displaystyle\frac{2\alpha_4}{\alpha_1\alpha_2}=\lambda,$
i.e., the family $\ll^{0,6}_{(\mu,1,0,\lambda,-1)}$, $\mu=2,
\lambda\in \CC\setminus\{0\}$ is obtained.
\end{itemize}
\end{itemize}

\

\noindent\fbox{\bf $\alpha_2\neq 0$,\ $\alpha_3\neq 0$ }
\begin{itemize}
\item $\alpha_1^2-4\alpha_3\neq 0$,\
$\alpha_1\alpha_2-2\alpha_3\neq 0$.
Taking appropriate $A_4$ and $B_4$ we derive
$$\alpha'_1=0, \ \ \alpha'_2=1, \quad \alpha'_3=-\displaystyle\frac{(\alpha_1\alpha_2-2\alpha_3)^2}{\alpha_2^2(\alpha_1^2-4\alpha_3)}=\mu,
$$
$$\alpha'_4=-\displaystyle\frac{(\alpha_1\alpha_2-2\alpha_3)(\alpha_1\alpha_2-2\alpha_4)}{\alpha_2^2(\alpha_1^2-4\alpha_3)}=\lambda.$$
Hence, we obtain the family of algebras
$\ll^{0,7}_{(0,1,\mu,\lambda,-1)}$, where $\mu\in
\CC\setminus\{0\}, \ \lambda\in\CC.$

\

\item $\alpha_1^2-4\alpha_3\neq 0$,\
$\alpha_1\alpha_2-2\alpha_3=0.$

It yields
$$\small\begin{array}{rll}
\alpha'_3-\alpha'_4&=&\displaystyle\frac{(\alpha_3-\alpha_4)\alpha_2A_1B_4^2}{(A_1+\alpha_2A_4)(\alpha_2A_1^2+2\alpha_3
A_1A_4+\alpha_2\alpha_3 A_4^2)},\\
\\
2\alpha'_3\alpha'_4-\alpha'^{2}_2\alpha'_3-\alpha'^{2}_4&=&\displaystyle\frac{(2\alpha_3\alpha_4-\alpha_2^2\alpha_3-\alpha_4^2)\alpha_2^2A_1^2B_4^4}{(A_1+\alpha_2A_4)^2
(\alpha_2A_1^2+2\alpha_3A_1A_4+\alpha_2\alpha_3A_4^2)^2}.
\end{array}
$$

\

\begin{itemize}
\item[$\checkmark$]$\alpha_3-\alpha_4=0$.

\

Therefore,
$2\alpha_3\alpha_4-\alpha_2^2\alpha_3-\alpha_4^2\neq 0$ and taking
the suitable values of $A_4$ and $B_4$ we obtain $\alpha'_1=4,$ $
\alpha'_2=1, $ $\alpha'_3=2$, $\alpha'_4=2$. Thus, the algebra
$\ll^{0,8}_{(-2\lambda,1,-\lambda,2,-1)}$ with $\lambda=-2$ is
obtained.

\

\item[$\checkmark$] $\alpha_3-\alpha_4\neq 0$.

\

If \mbox{$2\alpha_3\alpha_4-\alpha_2^2\alpha_3-\alpha_4^2=0,$ then
$\alpha_4\neq 0,$
$\alpha_3=\displaystyle\frac{\alpha_4^2}{2\alpha_4-\alpha_2^2}$,
$\alpha_4\neq \displaystyle\frac{\alpha_2^2}{2}$.}
Choosing adequate values of $A_4$ and $B_4$ we obtain
$\alpha'_1=8/3$, $ \alpha'_2=1,$ $\alpha'_3=4/3 $, $\alpha'_4=2$,
i.e., we derive the algebra $\ll^{0,8}_{(-2\lambda,1,-\lambda,2,-1)}$
with $\lambda=-4/3.$

\

If
$2\alpha_3\alpha_4-\alpha_2^2\alpha_3-\alpha_4^2\neq 0,$ then as
before we deduce $\alpha'_1=2\alpha'_3,$ $\alpha'_2=1, $
$\alpha'_3=-\displaystyle\frac{(\alpha_3-\alpha_4)^2}{2\alpha_3\alpha_4-\alpha_2^2\alpha_3-\alpha_4^2}=\lambda$,
$\alpha'_4=0$ and the family
$\ll^{0,9}_{(2\lambda,1,\lambda,0,-1)}$ with $\lambda\in
\CC\setminus\{0,1\}$ is obtained.
\end{itemize}

\

\item $\alpha_1^2-4\alpha_3=0,\ \alpha_1\alpha_2-2\alpha_3\neq
0$.

Then, $\alpha_1\neq 2\alpha_2$,\quad
$\alpha'^2_{1}-4\alpha'_4=\displaystyle\frac{(\alpha_1^2-4\alpha_4)4A_1B_4^2}{(2A_1+\alpha_1
A_4)^2(A_1+\alpha_2A_4)}.$

\

\begin{itemize}
\item[$\checkmark$] $\alpha_1^2-4\alpha_4=0$.

Then, $\alpha_1\alpha_2-2\alpha_4\neq 0$ and from the above
we deduce $\alpha'_1=1, \alpha'_2=1, \alpha'_3=1/4,
\alpha'_4=1/4.$ So, we obtain the algebra
$\ll^{0,10}_{(\lambda,1,\lambda^{2}/4,\mu,-1)}$ with $\lambda=1,
\mu=1/4.$

\item[$\checkmark$] $\alpha_1^2-4\alpha_4\neq 0,$ $\alpha_1\alpha_2-2\alpha_4= 0$
$\Rightarrow \ \alpha'_1=1, $ $ \alpha'_2=1,$ $ \alpha'_3=1/4$,
$\alpha'_4=1/2,$ i.e., we obtain
$\ll^{0,10}_{(\lambda,1,\lambda^{2}/4,\mu,-1)}$ with $\lambda=1,
\mu=1/2.$

\item[$\checkmark$] $\alpha_1^2-4\alpha_4\neq 0,$ $\alpha_1\alpha_2-2\alpha_4\neq
0$ $\Rightarrow \ \alpha'_1=1,$
$\alpha'_2=\displaystyle\frac{\alpha_1\alpha_2-2\alpha_4}{\alpha_1^2-4\alpha_4},$
$\alpha'_3=1/4$, $\alpha'_4=0$. The family
$\ll^{0,11}_{(1,\lambda,1/4,0,-1)}$, where $\lambda\in
\CC\setminus\{0,1/2\}$ is obtained.
\end{itemize}

\

\item $\alpha_1^2-4\alpha_3=0,\ \alpha_1\alpha_2-2\alpha_3= 0$.

Then, $\alpha_1=2\alpha_2,$\quad $\alpha_3=\alpha_2^2,$\quad
$\alpha'^2_{2}-\alpha'_4=\displaystyle\frac{(\alpha_2^2-\alpha_4)A_1B_4^2}{(A_1+\alpha_2A_4)^3}$.
\begin{itemize}
\item[$\checkmark$] $\alpha_2^2-\alpha_4= 0$.

Taking an appropriate value of $B_4$ it follows
that $\alpha'_1=2, $ $ \alpha'_2=1, $ $\alpha'_3=1,$ $
\alpha'_4=1$. Hence, we obtain
$\ll^{0,10}_{(\lambda,1,\lambda^{2}/4,\mu,-1)}$ with $\lambda=2,
\mu=1.$

\item[$\checkmark$] $\alpha_2^2-\alpha_4\neq 0$.

Choosing adequate $A_4$ and
$B_4=\displaystyle\frac{(\alpha_2^2-\alpha_4)A_1}{\alpha_2^3}$
yields $\alpha'_1=2,$ $ \alpha'_2=1,$ $\alpha'_3=1$ and
$\alpha'_4=0$. Thus, the algebra
$\ll^{0,10}_{(\lambda,1,\lambda^{2}/4,\mu,-1)}$ with $\lambda=2,
\mu=0$ is obtained.
\end{itemize}
\end{itemize}

Now, we consider the other case.

\

\fbox{{\bf Case 2.} $\epsilon=1$} ($n$ even)

\

Similar to the case 1, we apply the general change of generators
of basis. Then, we obtain all products and the following
expressions for $\alpha'_i, 1\leq i \leq 4$:
$$\small\begin{array}{l}
\alpha'_1=\displaystyle\frac{(A_1-A_4)(\alpha_1A_1+2\alpha_3A_4)}{A_1^2+\alpha_1
A_1A_4+\alpha_3 A_4^2},\quad \ \ \
\alpha'_2=\displaystyle\frac{\alpha_2(A_1-A_4)}{A_1+\alpha_2A_4},\\[5mm]
\alpha'_3=\displaystyle\frac{\alpha_3(A_1-A_4)^2}{A_1^2+\alpha_1
A_1A_4+\alpha_3 A_4^2}, \hspace{1.35cm}
\alpha'_4=\displaystyle\frac{(A_1-A_4)^2(\alpha_4 A_1+\alpha_2
\alpha_3 A_4)}{(A_1+\alpha_2A_4)(A_1^2+\alpha_1 A_1A_4+\alpha_3
A_4^2)},
\end{array}$$
verifying the restriction $
A_1(A_1-A_4)(A_1+\alpha_2A_4)(A_1^2+\alpha_1 A_1A_4+\alpha_3
A_4^2)\neq 0.$

Note that for these parameters we have
$$\begin{array}{lll}
\alpha'^2_{1}-4\alpha'_3&=&\displaystyle\frac{(\alpha_1^2-4\alpha_3)A_1^2(A_1-A_4)^2}{(A_1^2+\alpha_1
A_1A_4+\alpha_3 A_4^2)^2}, \\[5mm]
\alpha'_1\alpha'_2-2\alpha'_3&=&-\displaystyle\frac{(\alpha_1\alpha_2-2\alpha_3)A_1(A_1-A_4)^2}{(A_1+\alpha_2A_4)(A_1^2+\alpha_1
A_1A_4+\alpha_3 A_4^2)}, \\[5mm]
\alpha'_1\alpha'_2-2\alpha'_4&=&\displaystyle\frac{(\alpha_1\alpha_2-2\alpha_4)A_1(A_1-A_4)^2}{(A_1+\alpha_2A_4)(A_1^2+\alpha_1
A_1A_4+\alpha_3 A_4^2)}, \\[5mm]
\alpha'_1+2\alpha'_3&=&\displaystyle\frac{(\alpha_1+2\alpha_3)(A_1-A_4)A_1}{A_1^2+\alpha_1
A_1A_4+\alpha_3 A_4^2}.
\end{array}$$

Consequently, the nullity of the expressions
$\alpha_1^2-4\alpha_3, \alpha_1\alpha_2-2\alpha_3,
\alpha_1\alpha_2-2\alpha_4,\ \alpha_1+2\alpha_3$ are invariants.

Applying arguments as in the case 1 for the following subcases:

\

\noindent \fbox{\bf $\alpha_2=0$\ $\alpha_3=0$ }\ , \ \fbox{\bf
$\alpha_2=0$,\ $\alpha_3\neq 0$ }\ , \ \fbox{\bf $\alpha_2\neq
0$\ $\alpha_3=0$ }\ , \ \fbox{\bf $\alpha_2\neq 0$,\
$\alpha_3\neq 0$}\

\

\noindent we obtain the rest algebras and families of the theorem. $\hfill\Box$

\

The next theorem completes the classification of naturally graded
Leibniz algebras with characteristic sequence $(n-3,3)$.
\begin{thm}\label{first}
Let $\ll$ be an $n$-dimensional naturally graded Leibniz algebra
of the first type ($n\geq 9$). Then it is isomorphic to one of
the following pairwise non-isomorphic algebras:
$$\small\begin{array}{ll}
\ll^{34}_{(0,\lambda,0)}:&\ll^{35}_{(\mu,\lambda,1)}:\\[2mm]
\left\{
\begin{array}{l}
[e_i,e_1]=e_{i+1},\ 1\leq i \leq n-1,\ i\neq n-3,\\{}
[e_{1},e_{n-2}]=\lambda e_{n-1},\\{} [e_2,e_{n-2}]=\lambda e_{n},\
\lambda\in \CC.
\end{array}\right.&
\left\{
\begin{array}{l}
[e_i,e_1]=e_{i+1},\ 1\leq i \leq n-1,\ i\neq n-3,\\{}
[e_{1},e_{n-2}]=\mu e_2+\lambda e_{n-1},\\{} [e_2,e_{n-2}]=\mu
e_3+\lambda e_{n},\ \lambda,\mu \in \{0,1\}\\{} [e_{i},e_{n-2}]=\mu e_{i+1},\ 3 \leq i \leq
n-4,\\{} [e_{i},e_{n-2}]= e_{i+1},\ n-2\leq i\leq n-1,\
.
\end{array}\right.
\end{array}$$
$$\small\begin{array}{ll}
\ll^{36}_{(1,\lambda,0)}:&\ll^{37}_{(1,\lambda,2)}:\\[2mm]
\left\{
\begin{array}{l}
[e_i,e_1]=e_{i+1},\ 1\leq i \leq n-1,\ i\neq n-3,\\{}
[e_{1},e_{n-2}]=e_2+\lambda e_{n-1},\\{} [e_2,e_{n-2}]=e_3+\lambda
e_{n},\ \lambda\in \{-1,0\}\\{} [e_{i},e_{n-2}]=e_{i+1},\ 3 \leq i \leq n-4.
\end{array}\right.&
\left\{
\begin{array}{l}
[e_i,e_1]=e_{i+1},\ 1\leq i \leq n-1,\ i\neq n-3,\\{}
[e_{1},e_{n-2}]=e_2+\lambda e_{n-1},\\{}
[e_2,e_{n-2}]=e_3+\lambda
e_{n},\ \lambda\in \CC\\{}
[e_{i},e_{n-2}]=e_{i+1},\ 3 \leq i \leq n-4,\\{}
[e_{i},e_{n-2}]= 2e_{i+1},\ n-2\leq i\leq n-1.
\end{array}\right.
\end{array}$$
$$\small\begin{array}{ll}
\small\begin{array}{ll}
\ll^{38}_{(0,0,\lambda)}:&\ll^{39}_{(0,1,\lambda)}:\\[2mm]
\left\{
\begin{array}{l}
[e_{i}, e_{1}] = e_{i+1}, \ 1\leq i \leq
n-1, \ i\neq n-3,\\ {}
[e_{1},e_{n-2}] = -e_{n-1}, \\ {}
[e_{2},e_{n-2}] = -(1+\lambda)e_{n}, \\ {} [e_{1}, e_{n-1}] =
\lambda e_{n}, \ \lambda\in\CC.
\end{array}\right.&
\left\{
\begin{array}{l}
[e_{i}, e_{1}] = e_{i+1}, \ 1\leq i \leq
n-1,\ i\neq n-3, \\ {}
[e_{1},e_{n-2}] = -e_{n-1}, \\ {}
[e_{2},e_{n-2}] = -(1+\lambda)e_{n}, \\ {}
[e_{n-1}, e_{n-2}] =-e_{n}, \\ {}
[e_{1},e_{n-1}] = \lambda e_{n}, \\ {}
[e_{n-2},e_{n-1}] = e_{n}, \ \lambda\in\{-1,0\}.
\end{array}\right.
\end{array}
\end{array}$$
$$\small\begin{array}{ll}
\small\begin{array}{ll}
\ll^{40}_{(1,\lambda,\mu)}:&\ll^{41}_{(1,-1,\lambda)}:\\[2mm]
\left\{
\begin{array}{l}
[e_{i}, e_{1}] = e_{i+1}, \ 1\leq i \leq n-1,\ i \neq n-3, \\ {}
[e_{1},e_{n-2}] = e_{2}-e_{n-1}, \\ {}
[e_{2},e_{n-2}] = e_{3}-(1+\mu)e_{n}, \\ {}
[e_{i}, e_{n-2}] = e_{i+1},\ 3\leq i \leq
n-4,\\ {} [e_{n-1}, e_{n-2}] = -\lambda e_{n}, \\ {} [e_{1},
e_{n-1}] = \mu e_{n}, \\ {}
[e_{n-2},e_{n-1}] = \lambda e_{n},\ \lambda\in\{0,1\},\ \mu\in\CC.
\end{array}\right.&
\left\{
\begin{array}{l}
[e_{i}, e_{1}] = e_{i+1}, \ 1\leq i \leq n-1,\ i\neq n-3, \\ {}
[e_{1},e_{n-2}] = e_{2}-e_{n-1}, \\ {}
[e_{2},e_{n-2}] = e_{3}-(1+\lambda)e_{n}, \\ {}
[e_{i}, e_{n-2}] = e_{i+1},\ 3\leq i \leq
n-4,\\ {} [e_{n-1}, e_{n-2}] = e_{n}, \\ {} [e_{1}, e_{n-1}] =
\lambda e_{n}, \\ {}
[e_{n-2},e_{n-1}] = -e_{n},\ \lambda\in\{-1,0\}.
\end{array}\right.
\end{array}
\end{array}$$
\end{thm}
\dem Let $\ll$ be a Leibniz algebra of the first type. Then we
have the following multiplication: $$
\small \begin{array}{l}
[e_{i},
e_{1}]=e_{i+1},\ 1\leq i \leq n-4, \quad [e_{n-3}, e_{1}]=0,\\{}
[e_{n-2}, e_{1}]=e_{n-1},\ [e_{n-1},e_{1}]=e_{n},\ [e_{n},
e_{1}]=0.
\end{array}$$

It is not difficult to verify that $$\ll_{1}=<e_{1}, e_{n-2}>,\
\ll_{2}=<e_{2}, e_{n-1}>,\ \ll_{3}=<e_{3}, e_{n}>,\
\ll_{i}=<e_{i}>, \ 4\leq i \leq n-3$$ and
$\{e_{2},e_{3},\dots,e_{n-3}\}\subseteq R(\ll)$. Therefore, to
define the multiplication in $\ll$ it suffice to study the
multiplication of the element $e_{n-2}$ from the right side.

Introduce denotations
$[e_{1}, e_{n-2}]=\alpha _{1} e_{2}+\alpha _{2} e_{n-1}$, \
$[e_{n-2}, e_{n-2}]= \beta _{1}e_{2}+\beta _{2}e_{n-1}$,
$[e_{2}, e_{n-2}]=\gamma _{1} e_{3}+\gamma _{2} e_{n}$,
\ $[e_{n-1}, e_{n-2}]= \delta_{1}e_{3}+\delta _{2}e_{n}$.

Then to verify the Leibniz identity $[e_i, [e_j, e_k]] = [[e_i,
e_j], e_k]-[[e_i, e_k], e_j]$ it suffice to consider

\hspace{1.5cm} $ j=1 \mbox{ and } k=n-2,n-1, n;\quad j=n-2 \mbox{
and } k=1, n-1, n;\quad$

\hspace{1.5cm} $ j=n-1 \mbox{ and } k=1, n-2, n;\quad j=n \mbox{
and } k=1, n-2,n-1. $

\

We consider several cases.

\

\fbox{$e_{n-2} \in R(\ll).$}

\

Then $\{e_{2},e_{3},\dots,e_{n}\}\subseteq R(\ll)$ and,
consequently, we have $\alpha_k=$ $\beta_k=$ $\gamma_k=$ $\delta_k=0, \ 1
\leq k \leq 2 $. Thus, we obtain the algebra $\ll^{34}_{(0,0,0)}$.

\

\fbox{$e_{n-2} \notin R(\ll), \ e_{n-1} \in R(\ll).$}

\

Then, $e_n\in R(\ll)$ and 
$\gamma_k=\alpha_k, \ \delta_k=\beta_k,
\ 1\leq k \leq 2, \ \beta _{1}=0$.
Thus, the multiplication table of $\ll$ can be expressed in the form:
$$\begin{array}{llll} [e_{i}, e_{1}]& = &e_{i+1},& 1\leq i \leq
n-1,\ i\neq n-3, \\ {} [e_{1},e_{n-2}]& = &\alpha_{1}e_{2}+\alpha
_{2}e_{n-1},& \\{}
[e_{2},e_{n-2}]& = &\alpha_{1}e_{3}+\alpha_{2}e_{n},& \\{}
[e_{i}, e_{n-2}]& = &\alpha _{1}e_{i+1},& 3\leq i \leq
n-4,\\{}
[e_{i}, e_{n-2}]& = &\beta_2e_{i+1}, & n-2 \leq i \leq n-1.
\end{array}$$

Taking the general change of generators of basis:
$$e'_1=\displaystyle\sum_{i=1}^{n} A_i e_i,\ \
e'_{n-2}=\displaystyle\sum_{i=1}^{n} B_i e_i$$
we obtain the new basis
$\{e'_1,e'_2,\dots,e'_{n-1},e'_n\}$.

We compute all products and the new parameters are the following:
$$\small\begin{array}{lll}
\alpha'_1=\displaystyle\frac{\alpha_1B_{n-2}}{A_1+\alpha_1
A_{n-2}}, & \alpha'_2=\displaystyle \frac{A_1(\alpha_2 A_1+\beta_2
A_{n-2}-\alpha_1 A_{n-2} )}{(A_1+\alpha_1 A_{n-2})(A_1+ \beta_2
A_{n-2})},&
\beta'_2=\displaystyle\frac{\beta_2B_{n-2}}{A_1+\beta_2 A_{n-2}},
\end{array}$$

satisfying the restrictions
$$\begin{array}{ll}
A_1B_{n-2}(A_1+\alpha_1 A_{n-2})(A_1+\beta_2 A_{n-2})\neq 0,&\\
B_i=0,& 1 \leq i \leq n-6,\\
(\alpha_1-\beta_2)B_{i}=0, & n-5 \leq i \leq n-4,\\
\alpha_2B_{i}=0, & n-5 \leq i \leq n-4,
\end{array}$$
and $B_{n-2}(\beta_2A_{n-1}+\alpha_2A_2-\alpha_1A_{n-1})=B_{n-1}(\beta_2A_{n-2}+\alpha_2A_1-\alpha_1A_{n-2}).$

Note that only coefficients $A_1, A_{n-2}$, $B_{n-2}$ participate
in the expressions for the parameters $\alpha'_1$, $\alpha'_2$,
$\beta'_2$. Hence, we can suppose that $A_i=0, \ i\neq \{1,n-2\}$
and $B_j=0,$ with $j\neq n-2$.

It can be proved that the nullity of $\alpha_1-\beta_2$ is
invariant.

If $\alpha_1=0$, then the nullity of $1-\alpha_2$ is invariant and
if $\beta_2=0$, then the nullity of $1+\alpha_2$ is invariant, as
well.

Similar as in the proof of Theorem 2.1 we consider the
possible cases and in each of them we have the following pairwise
non-isomorphic algebras of the theorem:
$$\ll^{34}_{(0,\lambda,0)}, \ \lambda\in\CC-\{0\}; \ \
\ll^{35}_{(\mu,\lambda,1)}, \ \lambda,\mu\in\{0,1\};\ \
\ll^{36}_{(1,\lambda,0)},\ \lambda\in\{-1,0\};\ \
\ll^{37}_{(1,\lambda,2)},\ \lambda\in \CC. $$

\fbox{$e_{n-1} \notin R(\ll), \ e_{n} \in R(\ll).$}

\

Then, $e_{n-2}\notin R(\ll)$. Therefore, for defining the
multiplication of $\ll_1$ and $\ll_2$ it is enough to study the
multiplication of $e_{n-2}$ and $ e_{n-1}$ on the right side.

Introduce the notations
$$\begin{array}{ll}
[e_{1}, e_{n-2}]=\alpha _{1} e_{2}+\alpha _{2}
e_{n-1},&[e_{n-2}, e_{n-2}]= \beta _{1}e_{2}+\beta
_{2}e_{n-1},\\{}
[e_{2}, e_{n-2}]=\gamma _{1} e_{3}+\gamma_{2}e_{n},&
[e_{n-1}, e_{n-2}]=\delta_{1}e_{3}+\delta_{2}e_{n},\\{}
[e_{1}, e_{n-1}]=a _{1} e_{3}+a _{2} e_{n},&
[e_{n-2}, e_{n-1}]= b _{1}e_{3}+b _{2}e_{n}
\end{array}$$

From equality $[e_{i+1},e_{n-1}]=[[e_{i}, e_{n-1}], e_{1}], \
1\leq i \leq n-1 $ we obtain
$$\begin{array}{ll}
[e_{i},e_{n-1}] =a_1e_{i+2}, & 2 \leq i \leq n-5, \\ {}
[e_{i},e_{n-1}] = 0, & n-4 \leq i \leq n-3, \\ {}
[e_{n-1},e_{n-1}] = b_{1}e_{4}, & \\ {}
[e_{n},e_{n-1}] = b_{1}e_{5}. &
\end{array} $$

The equality $[e_{i+1},e_{n-2}]=[[e_{i}, e_{n-2}],
e_{1}]-[e_{i},e_{n-1}], \ 1\leq i \leq n-1$ yields
$$ \alpha_{2} = -1, \ \gamma_{1} = \alpha_{1}-a_{1}, \
\gamma_{2} = -(1+a_{2}),\
\delta_{1} = \beta_1-b_{1}, \
\delta_{2} = \beta_{2}-b_{2}, $$
$$[e_{i},e_{n-2}] = (\alpha_1-(i-1)a_{1})e_{i+1}, \quad 3 \leq i \leq n-4, $$
$$[e_{n-3},e_{n-2}] = 0, \
[e_{n},e_{n-2}] = (\beta_1-2b_{1})e_{4}.
$$

Since $[e_{1},e_{n-1}]\in R(\ll),$ then
$[[e_i,e_{n-1}],e_1]=[e_{i+1},e_{n-1}]$, which implies $b_1=0$
for $n\geq 9$. From $[e_n,[e_1,e_{n-2}]]=-[[e_n,e_{n-2}],e_1]$,
we obtain $\beta_1=0$. Consequently,  $[e_{n-2},e_{n-2}]\in
R(\ll)$ and $e_{n-1}\notin R(\ll)$ we have $\beta_2=0$.

Moreover, $[e_{n-2},e_{n-1}]\in R(\ll)$. Then
$[[e_{i},e_{n-2}],e_{n-1}]=[[e_{i},e_{n-1}],e_{n-2}],$ with $1 \leq i
\leq n$ and hence, $a_1=0$ for $n\geq 9$.

Thus, the multiplication in $\ll$ is as follows
$$\begin{array}{llll} [e_{i}, e_{1}] & = & e_{i+1}, & 1\leq i \leq
n-1,\ i\neq n-3, \\ {}
[e_{1},e_{n-2}]& = & \alpha_{1}e_{2}-e_{n-1},& \\ {}
[e_{2},e_{n-2}]& = & \alpha_{1}e_{3}-(1+a_{2})e_{n},& \\ {}
[e_{i}, e_{n-2}]& = & \alpha _{1}e_{i+1},& 3\leq i \leq
n-4,\\ {}
[e_{n-1}, e_{n-2}]& = &-b_{2}e_{n}, & \\ {}
[e_{1}, e_{n-1}]& = &a_{2}e_{n},& \\ {}
[e_{n-2},e_{n-1}] & = & b_{2}e_{n}. &
\end{array}$$

Similar as above we take the general change of generators of basis
and then we generate the new basis. After that we determine all products and the new
parameters are the following: $$\small\begin{array}{lll}
\alpha'_1=\displaystyle\frac{\alpha_1B_{n-2}}{A_1+\alpha_1
A_{n-2}}, &
b'_2=\displaystyle\frac{b_{2}B_{n-2}}{A_{1}-b_{2}A_{n-2}},&
a'_{2}=\displaystyle\frac{(a_{2}A_{1}+b_{2}A_{n-2})}{A_1-b_{2}A_{n-2}},
\end{array}$$
with the restrictions $$\begin{array}{ll}
&A_1B_{n-2}(A_1+\alpha_1 A_{n-2})(A_1-b_{2}A_{n-2})\neq 0,\\
&B_i=0,\qquad \quad 1 \leq i \leq n-3,\\ &
(1+a_{2})(-(A_{2}+\alpha_{1}A_{n-1})B_{n-2}+(A_1+\alpha_{1}A_{n-2})B_{n-1})=0.
\end{array}$$

Note that only coefficients $A_1, A_{n-2}$, $B_{n-2}$ participate
in the expressions for the parameters $\alpha'_1$, $b'_2$, $a'_2$.
Therefore, we can assume that $A_i=0,$ $i\neq 1,n-2$ and
$B_j=0,$ $j\neq n-2$.

It is proved that the nullity of $1+a_2,\ \alpha_1+b_2 $ are
invariants.
Similarly as in the proof of Theorem 2.1 we consider the
possible cases and for each of them we have the following pairwise
non-isomorphic Leibniz algebras:
$$\ll^{38}_{(0,0,\lambda)},\ \lambda\in\CC;\quad
\ll^{39}_{(0,1,\lambda)},\ \lambda\in\{-1,0\};$$
$$\ll^{40}_{(1,\lambda,\mu)}, \ \lambda\in\{0,1\},\ \mu\in\CC;\quad \ll^{41}_{(1,-1,\lambda)},\ \lambda\in\{-1,0\}. $$
\
\fbox{$e_{n }\notin R(\ll)$.}
\

Then, $e_{n-2}, \ e_{n-1}\notin R(\ll)$, as well.
We set \mbox{$[e_1,e_n]=\lambda_1
e_4, \ [e_{n-2},e_n]=\lambda_2 e_4.$}

From $[[e_i,e_1],e_n]=[[e_i,e_n],e_1], \ 1 \leq i \leq n$ it
follows that
$$\begin{array}{llll}
[e_i,e_n]&=&\lambda_1 e_{i+3}, & 2 \leq i \leq n-6, \\ {}
[e_i,e_n]&=&0, & n-5 \leq i \leq n-3, \\ {}
[e_{n-1},e_n]&=&\lambda_2
e_{5}, & \\ {} [e_{n},e_n]&=&\lambda_2 e_{6}. &
\end{array}
$$

However, $[[e_n,e_n],e_1]=0$ implies $\lambda_2=0$ for $n \geq 9$
and hence, $\lambda_1 \neq 0$.

\

If we denote $[e_{1}, e_{n-1}]=a _{1} e_{3}+a _{2} e_{n}$ and
$[e_{n-2}, e_{n-1}]= b _{1}e_{3}+b _{2}e_{n}.$  Then due to
$[e_1,e_{n-1}]+[e_{n-1},e_1]\in R(\ll)$ we have $a_2=-1$.

From
$[[e_i,e_1],e_{n-1}]=[[e_i,e_{n-1}],e_1]-[e_i,e_n], \ 1 \leq i
\leq n$ it follows that $b_1=0$ for $n\geq 9$ and $$
\begin{array}{llll}
[e_i,e_{n-1}] & = & (a_1-(i-1)\lambda_1) e_{i+2}, & 2 \leq i \leq
n-5,\\ {}
[e_i,e_{n-1}]& = & 0,& n-4 \leq i \leq n-3, \\ {}
[e_{n-2},e_{n-1}]& = & b_2 e_{n}, & \\ {}
[e_{n-1},e_{n-1}]&=&0,&\\ {}
[e_{n},e_{n-1}]&=&0. &
\end{array}
$$

The equality $[[e_i,e_n],e_{n-1}=[[e_i,e_{n-1}],e_n],\ 1 \leq i \leq n$
implies $\lambda_1=0$, which is a contradiction with condition
$\lambda_1\neq0$. Consequently, in this case there does not appear
any naturally graded Leibniz algebra. $\hfill\Box$

\

Summarizing the results of the Theorems 2.1. and 2.2 we complete
the classification of naturally graded Leibniz algebras with the characteristic sequence $(n-3,3)$.

\textbf{Acknowledgments.} \emph{This work is
supported in part by the PAI, FQM143 of Junta de
Andaluc\'{\i}a (Spain). B.A. Omirov was supported by a grant of
NATO-Reintegration ref. CBP.EAP.RIG. 983169 and he would like to thank of
the Universidad de Sevilla for their hospitality.}

{\sc Jes\'{u}s M. Cabezas.} Dpto. de Matem\'{a}tica Aplicada. Universidad del Pa\'{\i}s Vasco.
c) Nieves Cano 12, 01002 Vitoria. \ (Spain), \ e-mail:
\emph{jm.cabezas@ehu.es}
\\

{\sc Luisa M. Camacho, Jos\'{e} R. G\'{o}mez.} Dpto. Matem\'{a}tica Aplicada
I. Universidad de Sevilla. Avda. Reina Mercedes, s/n. 41012
Sevilla. (Spain), e-mail: \emph{lcamacho@us.es},
\emph{jrgomez@us.es}
\\

{\sc Bakhrom A. Omirov.} Institute of Mathematics and Information Technologues, Uzbekistan
Academy of Science, F. Hodjaev str. 29, 100125, Tashkent
(Uzbekistan), e-mail: \emph{omirovb@mail.ru}

\begin{thebibliography}{12}
\bibitem{Omirov1}
{\rm Ayupov Sh.A., Omirov B.A.} {\it On some classes of nilpotent
Leibniz algebras}, (Russian) Sibirsk. Mat. Zh., vol. 42 (1), 2001,
p. 18--29; translation in Siberian Math. J., vol. 42 (1), 2001,
p. 15--24.
\bibitem{J.Lie.Theory}{\rm Cabezas J.M., Pastor E.} {\it Naturally graded p-filiform Lie algebras
in arbitrary finite dimension}, Journal of Lie Theory, vol. 15,
2005, p. 379--391.
\bibitem{Comm-2-fil} {\rm Camacho L.M., G\'{o}mez J.R., Gonz\'{a}lez A.J., Omirov B.
A.} {\it Naturally graded 2-filiform Leibniz algebras}, to appear in Comm. Algebra.
\bibitem{Comm-p-fil} {\rm Camacho L.M., G\'{o}mez J.R., Gonz\'{a}lez A.J., Omirov B.
A.} {\it The classification of naturally graded p-filiform Leibniz algebras}, to appear in Comm. Algebra.
\bibitem{JSC} {\rm Camacho L.M., G\'{o}mez J.R., Gonz\'{a}lez A.J., Omirov B.
A.} {\it Naturally graded quasi-filiform Leibniz algebras},
Journal of Symbolic Computation, 44(5), 2009, p. 527-539.
\bibitem{Dz1}
{\rm Dzhumadil'daev A. S.} {\it Cohomologies of colour Leibniz
algebras: pre-simplicial approach.} Lie theory and its
applications in physics, III (Clausthal, 1999), p. 124--136, World
Sci. Publ., River Edge, NJ, 2000.
\bibitem{Dz2}
{\rm Dzhumadil'daev A. S., Davydov A. A.,} {\it Factor-complex for Leibniz cohomology.} Special issue dedicated to Alexei Ivanovich Kostrikin. Comm. Algebra, vol. 29(9), 2001, p.
4197--4210.
\bibitem{Cohomology}
{\rm Feigin B.L. and Fuks D.B.} {\it Cohomology of Lie groups and Lie algebras,} Itogi Nauki Tekh. Ser. Sovrem. Probl. Mat.
Fundam. Napravleniya 21 (1988), p. 121-209, Zbl. 653.17008.
\bibitem{fialiowski} {\rm Fialowski A.,} {\it On the cohomology of infinite-dimensional nilpotent Lie algebras,} Adv. Math., 97
(1993), p. 267-277.
\bibitem{G-H}
{\rm Goze M., Hakimjanov Y. (Khakimdjanov)}, {\it Nilpotent Lie
algebras}, Kluwer Academics Publishers, 1996.
\bibitem{loday} {\rm Loday J.L.} {\it Une version non commutative des
alg\`{e}bres de Lie: les alg\`{e}bres de Leibniz.} Ens. Math. vol. 39,
1993, p. 269--293.
\bibitem{loday2} {\rm Loday J.L.} {\it Overview on Leibniz algebras, dialgebras and their homology}, Fields Ins. Commun., vol. 181, 1997, p. 91--102.
\bibitem{Ve}
{\rm Vergne M.} {\it Cohomologie des alg\`{e}bres de Lie nilpotentes.
Application \`{a} l'\'{e}tude de la variet\'{e} des alg\`{e}bres de Lie
nilpotentes}, Bull. Soc. Math. France, vol. 98, 1970, p. 81--116.
\end{thebibliography}
\end{document}